\input amstex
\documentstyle{amsppt}
\NoBlackBoxes
\mag1100
\hyphenation{french}
\let\aro=@ 
\newcount\numcount
\def\numerote{\global\advance\numcount by 1 \the\numcount}
\def\lastnum[#1]{{\advance \numcount by #1(\the\numcount)}}
\def\temp{1.34}%
\let\tempp=\relax
\expandafter\ifx\csname psboxversion\endcsname\relax
  \message{PSBOX(\temp) loading}%
\else
    \ifdim\temp cm>\psboxversion cm
      \message{PSBOX(\temp) loading}%
    \else
      \message{PSBOX(\psboxversion) is already loaded: I won't load
        PSBOX(\temp)!}%
      \let\temp=\psboxversion
      \let\tempp= 
    \fi
\fi
\tempp
\let\psboxversion=\temp
\catcode`\@=11
%
%
\def\psfortextures{
\def\PSspeci@l##1##2{%
\special{illustration ##1\space scaled ##2}%
}}%
\def\psfordvitops{
\def\PSspeci@l##1##2{%
\special{dvitops: import ##1\space \the\drawingwd \the\drawinght}%
}}%
\def\psfordvips{
\def\PSspeci@l##1##2{%
\d@my=0.1bp \d@mx=\drawingwd \divide\d@mx by\d@my
\includegraphics{##1\space}}}%
\def\psforoztex{
\def\PSspeci@l##1##2{%
\special{##1 \space
      ##2 1000 div dup scale
      \number-\psllx\space \number-\pslly\space translate
}}}%
\def\psfordvitps{
\def\psdimt@n@sp##1{\d@mx=##1\relax\edef\psn@sp{\number\d@mx}}
\def\PSspeci@l##1##2{%
\special{dvitps: Include0 "psfig.psr"}
\psdimt@n@sp{\drawingwd}
\special{dvitps: Literal "\psn@sp\space"}
\psdimt@n@sp{\drawinght}
\special{dvitps: Literal "\psn@sp\space"}
\psdimt@n@sp{\psllx bp}
\special{dvitps: Literal "\psn@sp\space"}
\psdimt@n@sp{\pslly bp}
\special{dvitps: Literal "\psn@sp\space"}
\psdimt@n@sp{\psurx bp}
\special{dvitps: Literal "\psn@sp\space"}
\psdimt@n@sp{\psury bp}
\special{dvitps: Literal "\psn@sp\space startTexFig\space"}
\special{dvitps: Include1 "##1"}
\special{dvitps: Literal "endTexFig\space"}
}}%
\def\psfordvialw{
\def\PSspeci@l##1##2{
\special{language "PostScript",
position = "bottom left",
literal "  \psllx\space \pslly\space translate
  ##2 1000 div dup scale
  -\psllx\space -\pslly\space translate",
include "##1"}
}}%
\def\psforptips{
\def\PSspeci@l##1##2{{
\d@mx=\psurx bp
\advance \d@mx by -\psllx bp
\divide \d@mx by 1000\multiply\d@mx by \xscale
\incm{\d@mx}
\let\tmpx\dimincm
\d@my=\psury bp
\advance \d@my by -\pslly bp
\divide \d@my by 1000\multiply\d@my by \xscale
\incm{\d@my}
\let\tmpy\dimincm
\d@mx=-\psllx bp
\divide \d@mx by 1000\multiply\d@mx by \xscale
\d@my=-\pslly bp
\divide \d@my by 1000\multiply\d@my by \xscale
\at(\d@mx;\d@my){\special{ps:##1 x=\tmpx, y=\tmpy}}
}}}%
\def\psonlyboxes{
\def\PSspeci@l##1##2{%
\at(0cm;0cm){\boxit{\vbox to\drawinght
  {\vss\hbox to\drawingwd{\at(0cm;0cm){\hbox{({\tt##1})}}\hss}}}}
}}%
\def\psloc@lerr#1{%
\let\savedPSspeci@l=\PSspeci@l%
\def\PSspeci@l##1##2{%
\at(0cm;0cm){\boxit{\vbox to\drawinght
  {\vss\hbox to\drawingwd{\at(0cm;0cm){\hbox{({\tt##1}) #1}}\hss}}}}
\let\PSspeci@l=\savedPSspeci@l
}}%
%
%
\newread\pst@mpin
\newdimen\drawinght\newdimen\drawingwd
\newdimen\psxoffset\newdimen\psyoffset
\newbox\drawingBox
\newcount\xscale \newcount\yscale \newdimen\pscm\pscm=1cm
\newdimen\d@mx \newdimen\d@my
\newdimen\pswdincr \newdimen\pshtincr
\let\ps@nnotation=\relax
{\catcode`\|=0 |catcode`|\=12 |catcode`|
|catcode`#=12 |catcode`*=14
|xdef|backslashother{\}*
|xdef|percentother{
|xdef|tildeother{~}*
|xdef|sharpother{#}*
}%
\def\R@moveMeaningHeader#1:->{}%
\def\uncatcode#1{%
\edef#1{\expandafter\R@moveMeaningHeader\meaning#1}}%
\def\execute#1{#1}
\def\psm@keother#1{\catcode`#112\relax}
\def\executeinspecs#1{%
\execute{\begingroup\let\do\psm@keother\dospecials\catcode`\^^M=9#1\endgroup}}%
\def\@mpty{}%
\def\matchexpin#1#2{
  \fi%
  \edef\tmpb{{#2}}%
  \expandafter\makem@tchtmp\tmpb%
  \edef\tmpa{#1}\edef\tmpb{#2}%
  \expandafter\expandafter\expandafter\m@tchtmp\expandafter\tmpa\tmpb\endm@tch%
  \if\match%
}%
\def\matchin#1#2{%
  \fi%
  \makem@tchtmp{#2}%
  \m@tchtmp#1#2\endm@tch%
  \if\match%
}%
\def\makem@tchtmp#1{\def\m@tchtmp##1#1##2\endm@tch{%
  \def\tmpa{##1}\def\tmpb{##2}\let\m@tchtmp=\relax%
  \ifx\tmpb\@mpty\def\match{YN}%
  \else\def\match{YY}\fi%
}}%
\def\incm#1{{\psxoffset=1cm\d@my=#1
 \d@mx=\d@my
  \divide\d@mx by \psxoffset
  \xdef\dimincm{\number\d@mx.}
  \advance\d@my by -\number\d@mx cm
  \multiply\d@my by 100
 \d@mx=\d@my
  \divide\d@mx by \psxoffset
  \edef\dimincm{\dimincm\number\d@mx}
  \advance\d@my by -\number\d@mx cm
  \multiply\d@my by 100
 \d@mx=\d@my
  \divide\d@mx by \psxoffset
  \xdef\dimincm{\dimincm\number\d@mx}
}}%
%
\newif\ifNotB@undingBox
\newhelp\PShelp{Proceed: you'll have a 5cm square blank box instead of
your graphics (Jean Orloff).}%
\def\s@tsize#1 #2 #3 #4\@ndsize{
  \def\psllx{#1}\def\pslly{#2}%
  \def\psurx{#3}\def\psury{#4}
  \ifx\psurx\@mpty\NotB@undingBoxtrue
  \else
    \drawinght=#4bp\advance\drawinght by-#2bp
    \drawingwd=#3bp\advance\drawingwd by-#1bp
  \fi
  }%
\def\sc@nBBline#1:#2\@ndBBline{\edef\p@rameter{#1}\edef\v@lue{#2}}%
\def\g@bblefirstblank#1#2:{\ifx#1 \else#1\fi#2}%
{\catcode`\%=12
\xdef\B@undingBox{
\def\ReadPSize#1{
 \readfilename#1\relax
 \let\PSfilename=\lastreadfilename
 \openin\pst@mpin=#1\relax
 \ifeof\pst@mpin \errhelp=\PShelp
   \errmessage{I haven't found your postscript file (\PSfilename)}%
   \psloc@lerr{was not found}%
   \s@tsize 0 0 142 142\@ndsize
   \closein\pst@mpin
 \else
   \if\matchexpin{\GlobalInputList}{, \lastreadfilename}%
   \else\xdef\GlobalInputList{\GlobalInputList, \lastreadfilename}%
     \immediate\write\psbj@inaux{\lastreadfilename,}%
   \fi%
   \loop
     \executeinspecs{\catcode`\ =10\global\read\pst@mpin to\n@xtline}%
     \ifeof\pst@mpin
       \errhelp=\PShelp
       \errmessage{(\PSfilename) is not an Encapsulated PostScript File:
           I could not find any \B@undingBox: line.}%
       \edef\v@lue{0 0 142 142:}%
       \psloc@lerr{is not an EPSFile}%
       \NotB@undingBoxfalse
     \else
       \expandafter\sc@nBBline\n@xtline:\@ndBBline
       \ifx\p@rameter\B@undingBox\NotB@undingBoxfalse
         \edef\t@mp{%
           \expandafter\g@bblefirstblank\v@lue\space\space\space}%
         \expandafter\s@tsize\t@mp\@ndsize
       \else\NotB@undingBoxtrue
       \fi
     \fi
   \ifNotB@undingBox\repeat
   \closein\pst@mpin
 \fi
\message{#1}%
}%
%
%
\def\psboxto(#1;#2)#3{\vbox{%
   \ReadPSize{#3}%
   \advance\pswdincr by \drawingwd
   \advance\pshtincr by \drawinght
   \divide\pswdincr by 1000
   \divide\pshtincr by 1000
   \d@mx=#1
   \ifdim\d@mx=0pt\xscale=1000
         \else \xscale=\d@mx \divide \xscale by \pswdincr\fi
   \d@my=#2
   \ifdim\d@my=0pt\yscale=1000
         \else \yscale=\d@my \divide \yscale by \pshtincr\fi
   \ifnum\yscale=1000
         \else\ifnum\xscale=1000\xscale=\yscale
                    \else\ifnum\yscale<\xscale\xscale=\yscale\fi
              \fi
   \fi
   \divide\drawingwd by1000 \multiply\drawingwd by\xscale
   \divide\drawinght by1000 \multiply\drawinght by\xscale
   \divide\psxoffset by1000 \multiply\psxoffset by\xscale
   \divide\psyoffset by1000 \multiply\psyoffset by\xscale
   \global\divide\pscm by 1000
   \global\multiply\pscm by\xscale
   \multiply\pswdincr by\xscale \multiply\pshtincr by\xscale
   \ifdim\d@mx=0pt\d@mx=\pswdincr\fi
   \ifdim\d@my=0pt\d@my=\pshtincr\fi
   \message{scaled \the\xscale}%
 \hbox to\d@mx{\hss\vbox to\d@my{\vss
   \global\setbox\drawingBox=\hbox to 0pt{\kern\psxoffset\vbox to 0pt{%
      \kern-\psyoffset
      \PSspeci@l{\PSfilename}{\the\xscale}%
      \vss}\hss\ps@nnotation}%
   \global\wd\drawingBox=\the\pswdincr
   \global\ht\drawingBox=\the\pshtincr
   \global\drawingwd=\pswdincr
   \global\drawinght=\pshtincr
   \baselineskip=0pt
   \copy\drawingBox
 \vss}\hss}%
  \global\psxoffset=0pt
  \global\psyoffset=0pt
  \global\pswdincr=0pt
  \global\pshtincr=0pt 
  \global\pscm=1cm 
}}%
%
%
\def\psboxscaled#1#2{\vbox{%
  \ReadPSize{#2}%
  \xscale=#1
  \message{scaled \the\xscale}%
  \divide\pswdincr by 1000 \multiply\pswdincr by \xscale
  \divide\pshtincr by 1000 \multiply\pshtincr by \xscale
  \divide\psxoffset by1000 \multiply\psxoffset by\xscale
  \divide\psyoffset by1000 \multiply\psyoffset by\xscale
  \divide\drawingwd by1000 \multiply\drawingwd by\xscale
  \divide\drawinght by1000 \multiply\drawinght by\xscale
  \global\divide\pscm by 1000
  \global\multiply\pscm by\xscale
  \global\setbox\drawingBox=\hbox to 0pt{\kern\psxoffset\vbox to 0pt{%
     \kern-\psyoffset
     \PSspeci@l{\PSfilename}{\the\xscale}%
     \vss}\hss\ps@nnotation}%
  \advance\pswdincr by \drawingwd
  \advance\pshtincr by \drawinght
  \global\wd\drawingBox=\the\pswdincr
  \global\ht\drawingBox=\the\pshtincr
  \global\drawingwd=\pswdincr
  \global\drawinght=\pshtincr
  \baselineskip=0pt
  \copy\drawingBox
  \global\psxoffset=0pt
  \global\psyoffset=0pt
  \global\pswdincr=0pt
  \global\pshtincr=0pt 
  \global\pscm=1cm
}}%
%
\def\psbox#1{\psboxscaled{1000}{#1}}%
\newif\ifn@teof\n@teoftrue
\newif\ifc@ntrolline
\newif\ifmatch
\newread\j@insplitin
\newwrite\j@insplitout
\newwrite\psbj@inaux
\immediate\openout\psbj@inaux=psbjoin.aux
\immediate\write\psbj@inaux{\string\joinfiles}%
\immediate\write\psbj@inaux{\jobname,}%
%
%
\def\toother#1{\ifcat\relax#1\else\expandafter%
  \toother@ux\meaning#1\endtoother@ux\fi}%
\def\toother@ux#1 #2#3\endtoother@ux{\def\tmp{#3}%
  \ifx\tmp\@mpty\def\tmp{#2}\let\next=\relax%
  \else\def\next{\toother@ux#2#3\endtoother@ux}\fi%
\next}%
%
%
\let\readfilenamehook=\relax
\def\re@d{\expandafter\re@daux}
\def\re@daux{\futurelet\nextchar\stopre@dtest}%
\def\re@dnext{\xdef\lastreadfilename{\lastreadfilename\nextchar}%
  \afterassignment\re@d\let\nextchar}%
\def\stopre@d{\egroup\readfilenamehook}%
\def\stopre@dtest{%
  \ifcat\nextchar\relax\let\nextread\stopre@d
  \else
    \ifcat\nextchar\space\def\nextread{%
      \afterassignment\stopre@d\chardef\nextchar=`}%
    \else\let\nextread=\re@dnext
      \toother\nextchar
      \edef\nextchar{\tmp}%
    \fi
  \fi\nextread}%
\def\readfilename{\bgroup%
  \let\\=\backslashother \let\%=\percentother \let\~=\tildeother
  \let\#=\sharpother \xdef\lastreadfilename{}%
  \re@d}%
%
%
\xdef\GlobalInputList{\jobname}%
\def\psnewinput{%
  \def\readfilenamehook{
    \if\matchexpin{\GlobalInputList}{, \lastreadfilename}%
    \else\xdef\GlobalInputList{\GlobalInputList, \lastreadfilename}%
      \immediate\write\psbj@inaux{\lastreadfilename,}%
    \fi%
    \ps@ldinput\lastreadfilename\relax%
    \let\readfilenamehook=\relax%
  }\readfilename%
}%
\expandafter\ifx\csname @@input\endcsname\relax    
  \immediate\let\ps@ldinput=\input\def\input{\psnewinput}%
\else
  \immediate\let\ps@ldinput=\@@input
  \def\@@input{\psnewinput}%
\fi%
\def\nowarnopenout{%
 \def\warnopenout##1##2{%
   \readfilename##2\relax
   \message{\lastreadfilename}%
   \immediate\openout##1=\lastreadfilename\relax}}%
\def\warnopenout#1#2{%
 \readfilename#2\relax
 \def\t@mp{TrashMe,psbjoin.aux,psbjoint.tex,}\uncatcode\t@mp
 \if\matchexpin{\t@mp}{\lastreadfilename,}%
 \else
   \immediate\openin\pst@mpin=\lastreadfilename\relax
   \ifeof\pst@mpin
     \else
     \errhelp{If the content of this file is so precious to you, abort (ie
press x or e) and rename it before retrying.}%
     \errmessage{I'm just about to replace your file named \lastreadfilename}%
   \fi
   \immediate\closein\pst@mpin
 \fi
 \message{\lastreadfilename}%
 \immediate\openout#1=\lastreadfilename\relax}%
{\catcode`\%=12\catcode`\*=14
\gdef\splitfile#1{*
 \readfilename#1\relax
 \immediate\openin\j@insplitin=\lastreadfilename\relax
 \ifeof\j@insplitin
   \message{! I couldn't find and split \lastreadfilename!}*
 \else
   \immediate\openout\j@insplitout=TrashMe
   \message{< Splitting \lastreadfilename\space into}*
   \loop
     \ifeof\j@insplitin
       \immediate\closein\j@insplitin\n@teoffalse
     \else
       \n@teoftrue
       \executeinspecs{\global\read\j@insplitin to\spl@tinline\expandafter
         \ch@ckbeginnewfile\spl@tinline
       \ifc@ntrolline
       \else
         \toks0=\expandafter{\spl@tinline}*
         \immediate\write\j@insplitout{\the\toks0}*
       \fi
     \fi
   \ifn@teof\repeat
   \immediate\closeout\j@insplitout
 \fi\message{>}*
}*
\gdef\ch@ckbeginnewfile#1
 \def\t@mp{#1}*
 \ifx\@mpty\t@mp
   \def\t@mp{#3}*
   \ifx\@mpty\t@mp
     \global\c@ntrollinefalse
   \else
     \immediate\closeout\j@insplitout
     \warnopenout\j@insplitout{#2}*
     \global\c@ntrollinetrue
   \fi
 \else
   \global\c@ntrollinefalse
 \fi}*
\gdef\joinfiles#1\into#2{*
 \message{< Joining following files into}*
 \warnopenout\j@insplitout{#2}*
 \message{:}*
 {*
 \edef\w@##1{\immediate\write\j@insplitout{##1}}*
\w@{
\w@{
\w@{
\w@{
\w@{
\w@{
\w@{
\w@{
\w@{
\w@{
\w@{\string\input\space psbox.tex}*
\w@{\string\splitfile{\string\jobname}}*
\w@{\string\let\string\autojoin=\string\relax}*
}*
 \expandafter\tre@tfilelist#1, \endtre@t
 \immediate\closeout\j@insplitout
 \message{>}*
}*
\gdef\tre@tfilelist#1, #2\endtre@t{*
 \readfilename#1\relax
 \ifx\@mpty\lastreadfilename
 \else
   \immediate\openin\j@insplitin=\lastreadfilename\relax
   \ifeof\j@insplitin
     \errmessage{I couldn't find file \lastreadfilename}*
   \else
     \message{\lastreadfilename}*
     \immediate\write\j@insplitout{
     \executeinspecs{\global\read\j@insplitin to\oldj@ininline}*
     \loop
       \ifeof\j@insplitin\immediate\closein\j@insplitin\n@teoffalse
       \else\n@teoftrue
         \executeinspecs{\global\read\j@insplitin to\j@ininline}*
         \toks0=\expandafter{\oldj@ininline}*
         \let\oldj@ininline=\j@ininline
         \immediate\write\j@insplitout{\the\toks0}*
       \fi
     \ifn@teof
     \repeat
   \immediate\closein\j@insplitin
   \fi
   \tre@tfilelist#2, \endtre@t
 \fi}*
}%
\def\autojoin{%
 \immediate\write\psbj@inaux{\string\into{psbjoint.tex}}%
 \immediate\closeout\psbj@inaux
 \expandafter\joinfiles\GlobalInputList\into{psbjoint.tex}%
}%
%
%
%
\def\centinsert#1{\midinsert\line{\hss#1\hss}\endinsert}%
\def\psannotate#1#2{\vbox{%
  \def\ps@nnotation{#2\global\let\ps@nnotation=\relax}#1}}%
\def\pscaption#1#2{\vbox{%
   \setbox\drawingBox=#1
   \copy\drawingBox
   \vskip\baselineskip
   \vbox{\hsize=\wd\drawingBox\setbox0=\hbox{#2}%
     \ifdim\wd0>\hsize
       \noindent\unhbox0\tolerance=5000
    \else\centerline{\box0}%
    \fi
}}}%
%
\def\at(#1;#2)#3{\setbox0=\hbox{#3}\ht0=0pt\dp0=0pt
  \rlap{\kern#1\vbox to0pt{\kern-#2\box0\vss}}}%
%
\newdimen\gridht \newdimen\gridwd
\def\gridfill(#1;#2){%
  \setbox0=\hbox to 1\pscm
  {\vrule height1\pscm width.4pt\leaders\hrule\hfill}%
  \gridht=#1
  \divide\gridht by \ht0
  \multiply\gridht by \ht0
  \gridwd=#2
  \divide\gridwd by \wd0
  \multiply\gridwd by \wd0
  \advance \gridwd by \wd0
  \vbox to \gridht{\leaders\hbox to\gridwd{\leaders\box0\hfill}\vfill}}%
%
\def\fillinggrid{\at(0cm;0cm){\vbox{%
  \gridfill(\drawinght;\drawingwd)}}}%
%
%
\def\textleftof#1:{%
  \setbox1=#1
  \setbox0=\vbox\bgroup
    \advance\hsize by -\wd1 \advance\hsize by -2em}%
\def\textrightof#1:{%
  \setbox0=#1
  \setbox1=\vbox\bgroup
    \advance\hsize by -\wd0 \advance\hsize by -2em}%
\def\endtext{%
  \egroup
  \hbox to \hsize{\valign{\vfil##\vfil\cr%
\box0\cr%
\noalign{\hss}\box1\cr}}}%
%
\def\frameit#1#2#3{\hbox{\vrule width#1\vbox{%
  \hrule height#1\vskip#2\hbox{\hskip#2\vbox{#3}\hskip#2}%
        \vskip#2\hrule height#1}\vrule width#1}}%
\def\boxit#1{\frameit{0.4pt}{0pt}{#1}}%
\catcode`\@=12 
%
 \psfordvips   
 
\def \scaledpicture #1 by #2 (#3 scaled #4) 
{\dimen0=#1 \dimen1=#2 \divide\dimen0 by 1000\multiply\dimen0 by #4 
\divide\dimen1 by 1000\multiply\dimen1 by #4 
$$\psboxto(\dimen0;\dimen1){#3.ps}$$} 
\topmatter
\title 
D\'eformations feuillet\'ees des vari\'et\'es de Hopf (Foliated Deformations of Hopf Manifolds)
\endtitle 
\rightheadtext{d\'eformations feuillet\'ees}
\author
Laurent Meersseman,
Marcel Nicolau, 
Alberto Verjovsky
\endauthor
\date 18 f\'evrier 2009\enddate 
\address 
{Laurent Meersseman}\hfill\hfill\linebreak
\indent{I.M.B.}\hfill\hfill\linebreak
\indent{Universit\'e de Bourgogne}\hfill\hfill\linebreak
\indent{B.P. 47870}\hfill\hfill\linebreak
\indent{21078 Dijon Cedex}\hfill\hfill\linebreak
\indent{France}\hfill\hfill
\endaddress
\email laurent.meersseman\@u-bourgogne.fr \endemail
\address
{Marcel Nicolau}\hfill\hfill\linebreak
\indent{Departament de Matematiques}\hfill\hfill\linebreak
\indent{Universitat Auton\`oma de Barcelona}\hfill\hfill\linebreak
\indent{E 08193 Bellaterra}\hfill\hfill\linebreak
\indent{Espagne}\hfill\hfill
\endaddress
\email nicolau\@ mat.uab.es \endemail
\address
{Alberto Verjovsky}\hfill\hfill\linebreak \indent{Instituto de
Matem\'aticas de la UNAM, Unidad Cuernavaca}\hfill\hfill\linebreak
\indent{Apartado Postal 273-3, Admon. de correos
No.3}\hfill\hfill\linebreak \indent{Cuernavaca, Morelos,
M\'exico}\hfill\hfill
\endaddress
\email alberto\@matcuer.unam.mx \endemail

  
\thanks
Ces r\'esultats entrent dans le cadre du projet COMPLEXE (ANR-08-JCJC-0130-01) du premier auteur. Cette recherche a pu aboutir gr\^ace aux financements suivants : bourse FABER de la r\'egion Bourgogne pour le premier auteur, projet MTM2008-02294 du Ministerio de Ciencia e Innovaci\'{o}n d'Espagne
pour le second et, pour le troisi\`eme, CONACyT proyecto U1 55084 et PAPIIT (Universidad
Nacional Aut\'onoma de M\'exico) \# IN102108.
\endthanks

\abstract
Dans cet article, nous nous int\'eressons \`a une classe d'exemples bien particuliers de feuilletages \`a feuilles complexes, dont le
type diff\'eomorphe est fix\'e. Ils poss\`edent une unique feuille compacte et toutes les feuilles non compactes viennent s'accumuler
sur la feuille compacte. Nous montrons que la structure complexe le long des feuilles non compactes est fix\'ee par la structure complexe de la feuille
compacte. Et inversement nous montrons que la structure complexe d'une feuille non compacte suffit \`a d\'eterminer la structure complexe des autres feuilles.
Nous nous servons de ces r\'esultats pour discuter des d\'eformations feuillet\'ees des vari\'et\'es de Hopf, analogue feuillet\'e de
la notion de grande d\'eformation.
\medskip
In this article, we focus on a very special class of foliations with complex leaves whose diffeomorphism type is fixed. They have a unique
compact leaf and the non-compact leaves all accumulate onto it. We show that the complex structure along the non-compact leaves
is fixed by the complex structure of the compact leaf. Reciprocally, we prove that the complex structure along a non-compact leaf determines the complex structure along the other leaves. We apply these results to the study of foliated deformations of Hopf manifolds, a
foliated analogue to the notion of deformation in the large.
\endabstract

\endtopmatter

\def\C{{\Bbb C}}
\def\R{{\Bbb R}}
\def\Z{{\Bbb Z}}
\def\Q{{\Bbb Q}}
\def\N{{\Bbb N}}
\document
\head
{\bf 0. Introduction}
\endhead

Soit $(X,\Cal F)$ un feuilletage par vari\'et\'es complexes de dimension sup\'erieure ou \'egale \`a trois. On suppose que $\Cal F$ v\'erifie la propri\'et\'e suivante.

\definition{Hypoth\`ese diff\'erentiable}

\noindent Il existe une feuille compacte $L$ diff\'eomorphe \`a $\Bbb S^{2n-1}\times\Bbb S^1$ d'holonomie contractante $C^{\infty}$-plate.

\enddefinition

Sous cette hypoth\`ese, un voisinage de $L$ dans $X$ est enti\`erement d\'etermin\'e \`a hom\'eomorphisme feuillet\'e pr\`es ; et 
d\'etermin\'e \`a diff\'eomorphisme feuillet\'e pr\`es par la classe de conjugaison d'un g\'en\'erateur du groupe d'holonomie
(cf \cite{C-LN, Chapter IV, Theorem 2}).

\medskip

Nous nous int\'eressons dans cet article aux structures complexes que l'on peut mettre sur ce voisinage de $L$, et en
particulier aux interactions entre la structure complexe de la feuille et la structure complexe des feuilles non compactes.
\medskip

Voici un exemple d'un tel feuilletage : on quotiente $W=\Bbb C^n\setminus\{0\}\times \Bbb R$ par l'action engendr\'ee par
$$
(z,t)\in W\longmapsto (g(z),h(t))\in W
$$
o\`u $g$ est une contraction holomorphe et o\`u $h$ est un diff\'eomorphisme de $\Bbb R$ fixant $0$ et v\'erifiant $\vert h(t)\vert < \vert t\vert$ pour $t\not =0$ (respectivement 
$\vert h(t)\vert > \vert t\vert$ pour $t\not =0$). Dans ce cas, le re\-v\^e\-te\-ment universel holomorphe de la feuille
compacte est $\Bbb C^n\setminus\{0\}$.
\proclaim{D\'efinition}
Nous appellerons {\bf vari\'et\'e de Hopf} une vari\'et\'e compacte complexe de dimension $n>2$ diff\'eomorphe \`a $\Bbb S^{2n-1}\times\Bbb S^1$ et
de rev\^etement universel holomorphe $\Bbb C^n\setminus\{0\}$.
\endproclaim

Il s'agit de l'analogue, en dimension plus grande, des surfaces de Hopf primaires.
\medskip

Les feuilles non compactes sont quant \`a elles
des copies de $\Bbb C^n \setminus \{0\}$. Suivant l'hypoth\`ese v\'erifi\'ee par $h$ (contractante ou dilatante), il existe un voisinage $V$ de $L$ dans $X$ tel que l'intersection $C$ de toute feuille non compacte avec $V$ est de l'un des
deux types suivants :
\medskip
\noindent (i) $\underline{\text{type }\infty}$ : les feuilles $C$ sont biholomorphes \`a $\Bbb C^n \setminus K$, o\`u $K$ est un compact de $\Bbb C^n$ contenant $0$. De plus, $C$ s'accumule sur $L$ lorsqu'on s'approche de l'infini.

\noindent (ii) $\underline{\text{type }0}$ : les feuilles $C$ sont biholomorphes \`a $K\setminus\{0\}$, o\`u $K$ est un compact de $\Bbb C^n$ contenant $0$. De plus, $C$ s'accumule sur $L$ lorsqu'on s'approche de $0$.
\medskip

Le but de cet article est de montrer que, si la feuille compacte de $(X,\Cal F)$ est une Hopf ou si une des feuilles non compactes
est de type $0$ ou de type $\infty$, alors $(X,\Cal F)$ est essentiellement l'exemple pr\'ec\'edent. Plus pr\'ecis\'ement, 
nous d\'emontrons les r\'esultats suivants.

\proclaim{Th\'eor\`eme 1}
Soit $(X,\Cal F)$ un feuilletage par vari\'et\'es complexes de dimension sup\'erieure ou \'egale \`a trois v\'erifiant les hypoth\`eses diff\'erentiables pr\'ec\'edentes.

Si, dans un voisinage de la feuille compacte $L$, une feuille est de type $0$, ou de type $\infty$,
alors $L$ est une vari\'et\'e de Hopf.
\endproclaim

\proclaim{Th\'eor\`eme 2 (r\'eciproque)} 
Soit $(X,\Cal F)$ un feuilletage par vari\'et\'es complexes de dimension sup\'erieure ou \'egale \`a trois v\'erifiant les hypoth\`eses diff\'erentiables pr\'ec\'e\-den\-tes.

Si la feuille compacte $L$ est une vari\'et\'e de Hopf, alors dans un voisinage de $L$, toutes les feuilles non compactes sont de type $0$, ou toutes les feuilles non compactes sont de type $\infty$.
\endproclaim

\proclaim{Corollaire 1}
Sous les hypoth\`eses du th\'eor\`eme 1, toutes les feuilles non compactes sont de type $0$, ou toutes les feuilles non compactes sont de type $\infty$.
\endproclaim

Dans ces conditions, nous dirons que $\Cal F$ {\it est de type $0$} ou que $\Cal F$ {\it est de type $\infty$}.

\proclaim{Th\'eor\`eme 3 (uniformisation)}
Sous les hypoth\`eses \'equivalentes des th\'eor\`emes 1 et 2, si de plus le feuilletage est de type infini, alors $V$ est CR-isomorphe au quotient d'un voisinage ouvert connexe de $\tilde L\simeq \Bbb C^n\setminus \{0\}$ 
dans $\Bbb C^n\setminus\{0\} \times \Bbb R$ par l'action engendr\'ee par un couple $(g,h)$ o\`u $g$ est une contraction holomorphe et o\`u $h$ est un diff\'eomorphisme de $\Bbb R$ fixant $0$ et v\'erifiant $\vert h(t)\vert > \vert t\vert$ pour $t\not =0$.
\endproclaim

\proclaim{Corollaire 2}
Soient $(X,\Cal F)$ et $(X',\Cal F')$ deux feuilletages par vari\'et\'es complexes de m\^eme dimension v\'erifiant l'hypoth\`ese
diff\'erentiable. On les suppose de plus tout deux de type infini. 

Alors $\Cal F$ et $\Cal F'$ sont CR-isomorphes au voisinage des feuilles compactes si et seulement si les feuilles compactes sont 
biholomorphes.
\endproclaim

Nous ne savons pas si le th\'eor\`eme 3 et le corollaire 2 restent vrais lorsque le(s) feuilletage(s) est (sont) de type $0$.
\medskip

Dans la derni\`ere section, nous appliquons les r\'esultats et les techniques de l'article aux grandes d\'eformations de vari\'et\'es
de Hopf ainsi qu'\`a une version feuillet\'ee de grandes d\'eformations.
\medskip

L'ingr\'edient essentiel du papier est le lemme de compactification uniforme de \cite{M-V}, que nous rappelons en Section 3.
 
\head
{\bf 1. Pr\'eliminaires diff\'erentiables}
\endhead

Le lemme suivant est une simple reformulation de l'hypoth\`ese diff\'erentiable.

\proclaim{Lemme}
Sous les hypoth\`eses diff\'erentiables pr\'ec\'edentes, un voisinage $V$ de $L$ dans $X$ s'identifie diff\'erentiablement \`a un ouvert de la suspension d'holonomie de $\tilde L_{diff}\times \Bbb R$ (o\`u $\tilde L$ est le rev\^etement
universel diff\'erentiable de $L$). 
\endproclaim

Autrement dit, soit $h$ le morphisme d'holonomie de $L$, qu'on verra comme une application lisse de $\Bbb R$ dans $\Bbb R$ avec comme seul point fixe $0$. Consid\'erons l'action de $\Bbb Z$ sur $\tilde L\times \Bbb R$ (suspension d'holonomie) donn\'ee 
par
$$ 
(z,t)\in\tilde L\times\Bbb R \longmapsto (p\cdot z, h^{-p}(t))\in \tilde L\times \Bbb R
$$
pour $p\in\Bbb Z$, la notation $p\cdot $ d\'esignant l'action holomorphe de $p\in\Bbb Z\simeq \pi_1(L)$ sur le rev\^etement universel holomorphe $\tilde L$. Elle d\'efinit par passage au quotient une
vari\'et\'e $X_h$ feuillet\'ee par un feuilletage $\Cal F_h$ \`a feuilles complexes de codimension un. Comme $h$ ne poss\`ede qu'un point fixe, \`a savoir $0$, il y a une unique feuille compacte biholomorphe \`a $L$ (par abus de
notation nous la noterons encore $L$) et toutes
les feuilles non compactes sont biholomorphes \`a $\tilde L$ et viennent s'accumuler sur $L$.
\medskip

Remarquons que l'exemple de l'introduction est la version CR de cette construction (il suffit de remplacer $h^{-1}$ par $h$ ; l'exposant 
$-1$ vient de la d\'efinition du morphisme d'holonomie cf \cite{C-LN}, mais ne joue aucun r\^ole ici). Le lemme dit 
que, pour toute structure complexe, un voisinage ouvert connexe de $L$ dans $X$ s'identifie {\it diff\'erentiablement} \`a un voisinage ouvert connexe de $L$ dans $X_h$. La figure ci-dessous illustre la situation.
\medskip
Le dessin de gauche repr\'esente $\tilde L\times\Bbb R$, celui de droite $X_h$. La feuille $L$ est en gras et le voisinage $V$ est compris entre les deux cercles en pointill\'e. Sur le dessin de gauche, le rev\^etement universel
$\tilde L$ de $L$ est en gras, et le relev\'e $\tilde V$ de $V$ est l'ouvert contenant $\tilde L$ d\'elimit\'e par les pointill\'es.
Chaque droite verticale repr\'esente une copie de $\Bbb R^{2n}$ ; comme $\tilde L$ est diff\'eomorphe \`a $\Bbb R^{2n}\setminus \{0\}$,
on enl\`eve au dessin une ligne horizontale de z\'eros, indiqu\'ee en pointill\'es.

\medskip
\hfil\scaledpicture 11.0in by 6.1in (feuilletage scaled 300) \hfil
\medskip

Remarquons que $L$ disconnecte $V$. Dans la suite, nous travaillerons essentiellement sur un mod\`ele \`a bord et ne consid\'ererons que $W=\tilde L\times [0,\infty[$ et la restriction de $\tilde V$ \`a $W$. Par abus de notation, on 
continuera \`a noter $X$, $V$ et $X_h$ les vari\'et\'es \`a bord issues de ce mod\`ele.
\medskip

Int\'eressons-nous maintenant \`a la structure complexe standard existant sur ce feuilletage (celle d\'ecrite dans l'introduction). 
Puisque nous nous restreignons \`a un mod\`ele \`a bord, elle est donn\'ee comme quotient de 
$W=\Bbb C^n\setminus\{0\}\times \Bbb R$ par l'action engendr\'ee par
$$
(z,t)\in W\longmapsto (g(z),h(t))\in W
$$
o\`u $g$ est une contraction holomorphe et o\`u $h$ est un diff\'eomorphisme de $\Bbb R$ fixant $0$ et v\'erifiant $\vert h(t)\vert < \vert t\vert$ pour $t\not =0$ (respectivement 
$\vert h(t)\vert > \vert t\vert$ pour $t\not =0$).
\medskip

Les feuilles $\Bbb C^n \setminus \{0\}$ de $\tilde W$ ont deux bouts holomorphiquement distincts : le bout $0$ strictement
pseudo-convexe, et le bout infini strictement pseudo-concave. Supposons $h$ contractante. On peut pr\'eciser le dessin pr\'ec\'edent.

\medskip
\hfil\scaledpicture 10.5in by 6.4in (feuilletagecontrac scaled 300) \hfil
\medskip

C'est le bout infini des feuilles non compactes qui s'accumule sur la feuille compacte.
\medskip

Si au contraire $h$ est dilatante, le dessin devient :

\medskip
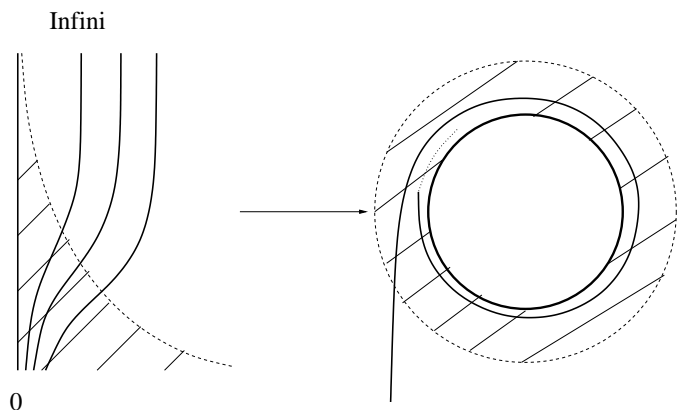
\hfil\scaledpicture 10.6in by 6.4in (feuilletagedil scaled 300) \hfil
\medskip

\noindent et c'est le bout $0$ des feuilles non compactes qui s'accumule sur la feuille compacte.
 
\head
{\bf 2. D\'eformations de vari\'et\'es complexes non compactes}
\endhead

Soit I un intervalle connexe de $\Bbb R$ contenant $0$ et soit $\pi : \Cal W\to I$ une famille de d\'eformations d'une vari\'et\'e complexe non compacte $W_0$. 

\remark{Remarque} 
Une telle famille sera toujours suppos\'ee diff\'erentiablement triviale, i.e. $\Cal W$ est diff\'eomorphe \`a $W_0\times I$.
\endremark
\medskip
Suivant la terminologie de \cite{A-V}, nous dirons que $\pi$ est {\it localement pseudo-triviale} si, \'etant donn\'e $K$ ouvert relativement compact de $W_0$, 
il existe $\Cal K$ ouvert relativement compact de $\Cal W$ tel que 
\medskip
\noindent (i) L'intersection de $\Cal K$ avec $W_0$ est $K$.

\noindent (ii) Il existe un CR-isomorphisme entre $\Cal K$ et $K\times J$ pour $J$ un voisinage ouvert connexe de $0$ dans $I$. 
\medskip

On peut dans ce contexte non compact prouver :

\proclaim{Proposition 1}
Soit $\pi : \Cal W\to I$ une famille de d\'eformations d'une vari\'et\'e complexe non compacte $W_0$. Soit $\Theta$ le faisceau des germes de champs tangents holomorphes sur $W_0$. Si le groupe de cohomologie
$H^1(W_0,\Theta)$ est nul, alors 
la famille $\pi$ est localement pseudo-triviale.
\endproclaim

Dans le cas $X$ compacte, l'annulation de $H^1(X,\Theta)$ entra\^{\i}ne la trivialit\'e locale de la famille de d\'eformations. Bien entendu, dans le cas non compact,
les fibres proches $W_t$ n'ont aucune raison d'\^etre biholomorphes \`a $W_0$ (penser aux exemples o\`u $W_0=\Bbb C$ et o\`u les autres fibres sont toutes des disques).

\demo{Preuve}
C'est une adaptation directe de la preuve de \cite{M-K, Theorem 3.2, p.45--55} pour le cas compact. Elle se fait en deux temps. Etant fix\'e un recouvrement local d'un voisinage $\Cal V$ de $W_0$  
par des cartes de submersion de $\Cal W$, on construit dans chaque carte une s\'erie formelle de telle sorte que ces s\'eries se recollent en un CR-isomorphisme formel entre $\Cal V$ et une d\'eformation triviale.
Ensuite, par des proc\'ed\'es de s\'eries majorantes, on montre la convergence de ce biholomorphisme formel quitte \`a restreindre $\Cal V$. 
\medskip
Dans le cas non compact, la partie formelle est exactement identique (la dif\-f\'e\-ren\-ce \'etant que le voisinage $\Cal V$ ne va a priori contenir qu'une unique fibre compl\`ete, \`a savoir $W_0$). Par contre,   
pour la preuve de la convergence, les majorations type (20) ainsi que les bornes uniformes pour certains cocycles donn\'ees par le lemme 3.7 de \cite{M-K} ne sont pas vraies sur un voisinage de $W_0$ tout entier, mais uniquement sur
un voisinage d'un compact de $W_0$.
$\square$
\enddemo

En particulier, on a

\proclaim{Corollaire 3}
Toute d\'eformation de $\Bbb C^n\setminus \{0\}$ v\'erifie les hypoth\`eses de la proposition 1 pour $n\geq 3$.
\endproclaim

\demo{Preuve}
Par \cite{Sc}, on a un isomorphisme entre $H^1(\Bbb C^n\setminus\{0\},\Theta)$ et $H^1(\Bbb C^n,\Theta)$ pour $n\geq 3$. Mais comme $\Bbb C^n$ est Stein, ce dernier groupe est nul.
$\square$
\enddemo

\head
{\bf 3. Feuilletages par vari\'et\'es complexes et lemme de compactification uniforme}
\endhead

Cette section est un rappel de d\'efinitions et de r\'esultats de \cite{M-V}.
\medskip

Un feuilletage par vari\'et\'es complexes $\Cal F$ de 
codimension un sur une vari\'et\'e lisse $M$ de dimension $2n+1$ est donn\'e par un atlas
$$ 
\Cal A=\{(U_i,\phi_i)_{i\in I}\quad\vert\quad \phi_i (U_i) \subset 
\Bbb C^n\times \Bbb R\cong \Bbb R^{2n+1}\} 
$$ 
tel que les changements de cartes
$$ 
(z,t)\in\phi_i(U_i\cap U_j)\longmapsto \phi_j 
\circ\phi_i^{-1}(z,t):=(\xi_{ij}(z,t),\zeta_{ij}(t))\in\phi_j(U_i\cap 
U_j) 
$$ 
sont holomorphes dans la direction tangente, {\it i.e.} l'application 
$\xi_{ij}$ est holomorphe \`a $t$ fix\'e. 
\medskip 
 Les plaques $\{t=\text{Cste}\}$ se recollent en des vari\'et\'es complexes connexes, les feuilles de $\Cal F$, qui forment
une partition de $M$.
\medskip
La structure complexe le long des feuilles est transversalement diff\'erentiable dans le sens suivant. L'op\'erateur presque complexe
$J$ d\'efini sur le fibr\'e tangent au feuilletage $T\Cal F$ est lisse non seulement le long des feuilles mais aussi transversalement.
\medskip
De fa\c con \'equivalente, un feuilletage par vari\'et\'es complexes sur $M$ peut aussi \^etre d\'efini \`a partir de la donn\'ee 
d'un sous-fibr\'e lisse $H$ du fibr\'e tangent $TM$ et d'une structure presque CR lisse $J$ sur $H$ tels que $H$ soit int\'egrable
(au sens de Frobenius, i.e. tangent \`a un feuilletage lisse) et $J$ soit int\'egrable (au sens complexe) le long des feuilles ;
ou encore tels que
$J$ soit int\'egrable et Levi-plate le long des feuilles.
\medskip

 
 

Venons-en au lemme de compactification uniforme de \cite{M-V}. Nous nous contenterons de l'\'enoncer dans le cas particulier d'un 
feuilletage v\'erifiant l'hypoth\`ese diff\'erentiable pr\'ecis\'ee dans l'introduction. Rappelons tout d'abord que si $L$ est une
vari\'et\'e complexe non compacte de dimension $n$, si $E$ est un bout de $L$ et $H$ une vari\'et\'e complexe de dimension strictement
inf\'erieure \`a $n$, alors une $E$-compactification de $L$ par $H$ est la donn\'ee d'une structure complexe sur $L\sqcup H$ telle que
\medskip
\noindent (i) Les inclusions naturelles de $L$ et $H$ dans $L\sqcup H$ soient holomorphes pour cette structure.

\noindent (ii) $H$ soit l'ensemble limite de $E$.
 \medskip

\proclaim{Lemme de compactification uniforme}
Soit $(X,\Cal F)$ un feuilletage par vari\'et\'es complexes de dimension sup\'erieure ou \'egale \`a trois v\'erifiant l'hypoth\`ese
dif\-f\'e\-ren\-ti\-a\-ble. 
Soit $V$ un voisinage de la feuille compacte $L$ pour lequel la conclusion du lemme de la section 1 s'applique. 
Soit $E$ un bout du rev\^etement universel $\tilde L$ 
de la feuille compacte, et par abus de notation le bout correspondant des feuilles non compactes dans $\tilde V$.
\medskip
Si $\tilde L$ {\bf ou} si une feuille non compacte de $V$
admet une $E$-compactification par une vari\'et\'e $H$, alors $\tilde L$ {\bf et toutes} les feuilles non compactes de $\tilde V$ 
admettent une $E$-compactification par $H$.
\medskip
De surcro\^{\i}t, cette compactification est uniforme au sens suivant : il existe une structure CR lisse sur 
$\tilde V\sqcup (H\times\Bbb R)$
telle que les injections naturelles de $\tilde V$ et $H\times\Bbb R$ dans $\tilde V\sqcup (H\times\Bbb R)$ soient CR et lisses.
\endproclaim  
 
\demo{Preuve}
C'est une application directe de la Proposition 5 de \cite{M-V}, apr\`es avoir not\'e que le feuilletage peut \^etre suppos\'e
``tame'' par \cite{M-V, Corollary 3} (la d\'efinition de tame dans ce contexte est donn\'ee dans \cite{M-V}). 
La seule diff\'erence est que nous permettons des compactifications  par des sous-vari\'et\'es
de codimension arbitraire, alors que tous les r\'esultats de compactification de \cite{M-V} sont faits dans le cadre des hypersurfaces.
Cependant, les preuves n'utilisent \`a aucun moment l'hypoth\`ese de codimension un et fonctionnent verbatim en toute g\'en\'eralit\'e.
$\square$
\enddemo

\head 
{\bf 4. Preuve du th\'eor\`eme 1}
\endhead

On consid\`ere le mod\`ele \`a bord de la section 1. 
Soit $N$ la trace d'une feuille non compacte de $\Cal F$ dans
$V$. Alors $N$ a deux bouts, de m\^eme que $\tilde L$, et l'un des bouts de $N$ s'accumule sur la feuille compacte $L$. Distinguons les cas (i) et (ii) du th\'eor\`eme.
\medskip
\noindent $\underline{\text{Type } 0}$ :
Par hypoth\`ese $N$ est de la forme $U\setminus \{0\}$ avec $U$ ouvert relativement compact de $\Bbb C^n$ contenant $0$. 
Et c'est le bout $0$ de $N$ qui s'accumule sur la feuille compacte $L$. 
\medskip
Ceci entra\^{\i}ne que $N$ peut \^etre compactifi\'ee holomorphiquement par un point le long du bout $0$. Le lemme de compactification 
implique qu'il en va de m\^eme pour $\tilde L$. 
Autrement dit, il existe une structure de vari\'et\'e complexe sur $\tilde L\sqcup\{0\}$ telle que l'injection naturelle $\tilde L\hookrightarrow \tilde L\sqcup\{0\}$ est holomorphe.
\medskip
Soit $g$ un g\'en\'erateur de l'action du groupe fondamental de $L$ sur $\tilde L$. Alors $g$ va se prolonger en un automorphisme de $\tilde L\sqcup\{0\}$ fixant $0$. L'application $g$ est contractante au sens de \cite{Ka}, puisqu'il s'agit d'une propri\'et\'e topologique
v\'erifi\'ee par le rev\^etement topologique de $L$.
Comme $\tilde L\sqcup\{0\}$ est lisse en $0$, toujours par \cite{Ka, Theorem 1}, on en d\'eduit que $\tilde L$ est $\Bbb C^n$, et donc que $L$ est une vari\'et\'e de Hopf.
\medskip
\noindent $\underline{\text{Type }\infty }$ :
Par hypoth\`ese, $N$ est de la forme $\Bbb C^n\setminus K$, pour $K$ un compact de $\Bbb C^n$. 
Et c'est le bout ``infini'' de $N$, not\'e $\infty$,
qui s'accumule sur la feuille $L$.
\medskip
Le lemme de compactification implique que $\tilde L$ peut \^etre compactifi\'ee \`a l'infini en ajoutant une copie de $\Bbb P^{n-1}$, i.e. il existe une structure de vari\'et\'e complexe
sur $\tilde L\sqcup \Bbb P^{n-1}$ telle que les injections naturelles $\tilde L\hookrightarrow \tilde L\sqcup \Bbb P^{n-1}$ et $\Bbb P^{n-1}\hookrightarrow \tilde L\sqcup \Bbb P^{n-1}$ soient holomorphes. Observons
que le fibr\'e normal de $\Bbb P^{n-1}$ dans $\tilde L\sqcup \Bbb P^{n-1}$ est le m\^eme que le fibr\'e normal de $\Bbb P^{n-1}$ dans $\Bbb P^n$, \`a savoir isomorphe \`a $\Cal O(1)$. En effet, toujours par le lemme de compactification, toutes les fibres de $\tilde V$ peuvent
\^etre compactifi\'ees par $\Bbb P^{n-1}$. L'ensemble des classes de Chern des fibr\'es normaux peut \^etre identifi\'e \`a une famille
de classes de cohomologie enti\`eres de $\Bbb P^{n-1}$ param\'etr\'ee par $t\in\Bbb R$. En r\'efl\'echissant un peu, on voit que la
propri\'et\'e d'uniformit\'e du lemme de compactification implique que cette famille est continue en $t$, donc constante de nombre de
Chern associ\'e constant \'egal \`a un. 

\medskip
De ce fait, un voisinage tubulaire $T$ bien choisi 
de $\Bbb P^{n-1}$ dans $\tilde L\sqcup \Bbb P^{n-1}$ a fronti\`ere $\partial T$ strictement pseudo-concave du c\^ot\'e de $\Bbb P^{n-1}$ d'apr\`es \cite{Ro, \S 5.2}, et donc strictement pseudo-convexe de l'autre c\^ot\'e. Toujours par
\cite{Ro, \S 4.3, 4.4}, il existe donc $W_0$ un espace analytique de Stein et une injection holomorphe $i$ d'un voisinage de cette fronti\`ere dans $W_0$. Ainsi en recollant le voisinage $T$ \`a l'ouvert  strictement pseudo-convexe de
$W_0$ bord\'e par l'image $i(\partial T)$ le long d'un voisinage de $\partial T$, on obtient un espace analytique compact $M$ contenant une sous-vari\'et\'e $H=\Bbb P^{n-1}$ \`a fibr\'e normal $\Cal O(1)$ tel que $M\setminus H$ soit Stein.
La figure ci-dessous illustre la construction de $M$.

\medskip

Remarquons que $M$ peut \^etre suppos\'ee plong\'ee dans $\Bbb P ^N$ pour $N$ grand avec $M\setminus H$ trace de ce plongement dans $\Bbb C^N$. En effet, il suffit de plonger un voisinage de $H$ dans $\Bbb P^N$ envoyant $H$ sur 
l'hyperplan \`a l'infini par \cite{Ro, Theorem 5.3} et d'\'etendre ce plongement \`a tout $M$ par \cite{G-R, Corollary VII.D.7}.
 
\medskip
\hfil\scaledpicture 8.6in by 7.8in (rec scaled 300) \hfil
\medskip

Nous affirmons qu'il existe une injection holomorphe de
$\tilde L\sqcup \Bbb P^{n-1}$ dans $M$ et que la diff\'erence entre ces deux espaces est constitu\'ee d'un unique point, appelons-le $z_0$. La preuve de cette affirmation s'inspire de \cite{Ka, p.572}.

\medskip

Pour simplifier, utilisons \'egalement la notation $i$ pour d\'efinir l'injection holomorphe de $T$ dans $M$. Soit $g$ la contraction holomorphe d\'efinie sur $\tilde L$. Il existe $T'\subset T$ voisinage strictement pseudo-concave de $\Bbb P^{n-1}$
tel que $g(T'\setminus\Bbb P^{n-1})$ soit contenu dans $T$. Via $i$, on transporte $g$ en un biholomorphisme $h$ d\'efini par
$$
h\equiv i\circ g\circ i^{-1}\ :\ i(T'\setminus\Bbb P^{n-1})\longrightarrow i(g(T'\setminus\Bbb P^{n-1}))
$$

Comme $i(T')$ est strictement pseudo-convexe, $h$ s'\'etend en un biholomorphisme de $M$ priv\'e de $i(\Bbb P^{n-1})$. On peut alors \'etendre l'injection holomorphe \`a 
$\tilde L$ tout entier en posant :

$$
i(z)=h^{p}\circ i\circ g^{-p}(z) \quad\text{pour}\quad z\in\tilde L,\ g^{-p}(z)\in T
$$
qui est bien d\'efinie gr\^ace \`a la relation
$$
i\circ g\equiv h\circ i \quad\text{sur}\quad g(T'\setminus\Bbb P^{n-1})
$$
et qui est d\'efini sur tout $\tilde L$ parce que $T'\setminus\Bbb P^{n-1}$ contient un domaine fondamental pour l'action induite par $g$. 
\medskip

Consid\'erons maintenant $M\setminus i(\Bbb P^{n-1})$ plong\'e dans $\Bbb C^N$ pour un grand $N$ et posons $S=M\setminus i(\tilde L\sqcup\Bbb P^{n-1})$. Nous allons montrer qu'il est r\'eduit \`a un point. C'est un compact dont on peut supposer qu'il contient 
l'origine de $\Bbb C^n$. En fait, la preuve aura pour cons\'equence que c'est l'ensemble des points fixes
de $h$. Notons que cela implique {\it a posteriori} qu'il s'agit  d'un ensemble analytique compact connexe dans $M\setminus i(\Bbb P^{n-1})$ Stein, donc d'un point. Cependant, nous allons devoir proc\'eder autrement.
\medskip
On remarque que :
\medskip
\noindent (i) Pour tout $z\in M\setminus S$, pour tout $p\in\Bbb Z$, on a $h^p(z)\in M\setminus S$.

\noindent (ii) Pour tout $z\in M\setminus S$, le point $z_{\infty}=\lim_{p\to\infty} h^p(z)$ qui est bien d\'efini par compacit\'e de $S$ (quitte \`a prendre une sous-suite) appartient \`a $S$ ; tandis que
$z_{-\infty}=\lim_{p\to -\infty} h^p(z)$ qui est bien d\'efini par compacit\'e de $i(\Bbb P^{n-1})$ (quitte \`a prendre une sous-suite) appartient \`a $\Bbb P^{n-1}$.
\medskip

On en d\'eduit que la fronti\`ere de $S$ est constitu\'e des limites quand $p$ va \`a l'infini des orbites $h^p(z)$. De surcro\^{\i}t, pour tout $p\in\Bbb Z$, on a $h^p(S)=S$. 
Par compacit\'e de $S$, quitte \`a
extraire des sous-suites, les suites $(h^p)_{p\geq 0}$ et $(h^p)_{p\leq 0}$ convergent uniform\'ement vers une fonction $h_0$ et son inverse d\'efinies sur $S$. D'autre part, comme $h$ est une contraction, la suite
 $(h^p)_{p\geq 0}$ converge uniform\'ement vers une fonction que nous continuerons d'appeler $h_0$ sur tout compact $S'$ de $M$ contenant $S$. Et d'apr\`es ce qui pr\'ec\`ede, $h_0$ envoie $S'$ sur $S$. Supposons de plus $S'$ connexe. 
Comme $h_0\equiv h_0\circ h^p$ pour tout $p>0$, $h_0$ atteint son maximum en un point int\'erieur \`a $S'$. Donc $h_0$ est constante sur $S'$ et aussi sur $S$. Mais $h_0$ est inversible sur $S$.
Autrement dit, $S$ est r\'eduit \`a un point.
\medskip

Ainsi $M$ est topologiquement $\Bbb P^n$ ; et en tant qu'espace complexe, $M$ est une vari\'et\'e projective qui poss\`ede 
(\'eventuellement) une singularit\'e en $z_0$. Soit $D$ le diviseur associ\'e \`a $\Bbb P^{n-1}$.

\proclaim{Lemme}
La vari\'et\'e projective $M$ v\'erifie

\noindent (i) Son anneau de cohomologie enti\`ere est isomorphe \`a 
$$
\Bbb Z\oplus\Bbb Z\cdot \omega\oplus\hdots\oplus\Bbb Z\cdot \omega^n
$$
pour $\omega$ un g\'en\'erateur du second groupe de cohomologie.

\noindent (ii) La section nulle du fibr\'e en droites associ\'e \`a $-D$ est n\'egative au sens de Grauert.

\noindent (iii) La classe de Chern de $D$ est $\omega$.

\noindent (iv) La dimension des sections globales des puissances positives de $D$ est donn\'ee par la formule
$$
\dim H^0(M,\Cal O(\nu D))=\pmatrix n+\nu \\ n \endpmatrix
$$
\endproclaim

\demo{Preuve}
La description topologique de la paire $(M,D)$ entra\^{\i}ne imm\'ediatement (i) et (iii) (quitte \`a remplacer $\omega$ par $-\omega$). Et nous avons d\'ej\`a vu que le fibr\'e normal de $D$ dans $M$ est $\Cal O(1)$, ce qui implique (ii) par
\cite{Ro, \S 5.2}. Finalement, on obtient (iv) par un calcul direct. Plus pr\'ecis\'ement, partant de la suite exacte courte
$$
0 \aro >>> \Cal O_M(-D) \aro >>> \Cal O_M \aro >>> \Cal O_M\vert _D \aro >>>0
$$
et tensorisant par $\Cal O_M(\nu D)$, on obtient la suite courte exacte
$$
 0 \aro >>> \Cal O_M((\nu -1) D) \aro >>> \Cal O_M (\nu D) \aro >>> \Cal O_D(\nu D) \aro >>>0
$$
Par simple connexit\'e de $M$, le premier groupe de cohomologie de $M$ \`a valeurs dans $\Cal O$ est nul. On d\'eduit de la suite exacte longue associ\'ee \`a la suite pr\'ec\'edente pour $\nu=1$ que
$$
H^1(M,\Cal O_M(D))=\{0\} \quad\text{ et }\quad H^0(M,\Cal O_M(D))=\Bbb C^{n+1}
$$
Par r\'ecurrence sur $\nu$, on compl\`ete le calcul des $H^0(M,\Cal O_M(\nu D))$.
$\square$
\enddemo

Pour conclure, on utilise le th\'eor\`eme 6 de \cite{H-K} qui affirme qu'une vari\'et\'e complexe v\'erifiant (i)-(ii)'-(iii)-(iv), avec
\medskip

\noindent (ii)' $M$ est K\"ahler de g\'en\'erateur $\omega$.
\medskip

\noindent est isomorphe \`a $\Bbb P^n$. On remarque en effet que, dans la preuve de \cite{H-K}, la condition (ii)' est uniquement utile pour plonger $M$ dans un projectif de grande dimension via le th\'eor\`eme de plongement de Kodaira
appliqu\'e \`a $\Cal O(D)$. Mais si l'on remplace le th\'eor\`eme de Kodaira par sa variante singuli\`ere de Grauert, la condition (ii) suffit.
\medskip

D\`es lors, $M$ est $\Bbb P^{n}$, donc $\tilde L$ est $\Bbb C^{n}\setminus\{0\}$ et $L$ est une Hopf.
$\square$

\head 
{\bf 5. Preuve des th\'eor\`emes 2 et 3}
\endhead
On se place dans un voisinage $V$ de la feuille compacte $L$ satisfaisant la conclusion du lemme de la section 1. 
Consid\'erons les feuilles non compactes. Elles s'accumulent sur la feuille compacte. Quand on passe au rev\^etement universel $\tilde V$, 
l'action du groupe fondamental envoie les
feuilles non compactes vers un bout du rev\^{e}tement universel $\Bbb C^{n}\setminus\{0\}$ de la feuille compacte, \`a savoir :
\medskip

\noindent (i) le bout $0$.

\noindent (ii) le bout infini.
\medskip

\noindent $\underline{\text{Cas (i)}}$ : le lemme de compactification uniforme implique que la famille des rev\^{e}\-te\-ments universels $\tilde V$ peut \^{e}tre
compactifi\'ee uniform\'ement en ajoutant un point \`a chaque fibre. On obtient ainsi une famille de d\'eformations $\Cal W$ dont la fibre centrale est $\tilde L\sqcup \{0\}\simeq \Bbb C^{n}$. 
D'apr\`es la proposition 1, elle est localement pseudo-triviale. En particulier, il existe $\Cal V$ voisinage de $0\in\tilde L\sqcup\{0\}$ dans 
$\Cal W$ et $\phi$ CR-isomorphisme de $\Cal V$ dans un voisinage de $\{0\}\times\{0\}$ dans $\tilde L\sqcup \{0\}\times [0,\infty [$. La compactification \'etant uniforme, $\phi$ envoie le voisinage \'epoint\'e 
de $0\in\tilde L\sqcup\{0\}$ sur un voisinage \'epoint\'e de $\{0\}\times\{0\}$ dans $\tilde L\times [0,\infty[$. Ce qui prouve le th\'eor\`eme 2 dans ce cas.
 
\medskip

\noindent $\underline{\text{Cas (ii)}}$ : le lemme de compactification uniforme implique que la famille $\tilde V$ des rev\^{e}tements universels peut \^{e}tre
compactifi\'ee uniform\'ement en ajoutant un diviseur isomorphe \`a $\Bbb P^{n-1}$ \`a chaque fibre. On obtient ainsi une famille $\Cal W$ de fibre centrale $\Bbb P^n\setminus\{0\}$.
\medskip

D'autre part, on a encore la pseudo-trivialit\'e locale de la famille $\tilde V$ d'apr\`es le corollaire 3.


\medskip

Ainsi, comme dans la section 4, on peut pratiquer la chirurgie suivante. On enl\`eve \`a $\Bbb P^n\setminus\{0\}$ (respectivement \`a $\Cal W$) un voisinage du bout $0$ de $\Bbb P^n\setminus \{0\}$ (respectivement du bout $0$ de $\Cal W$)
et on recolle un voisinage de $0$ dans $\Bbb C^n$ (respectivement d'un produit de ce voisinage pare $[0,\epsilon[$) le long d'un anneau (respectivement du produit d'un anneau par $[0,\epsilon[$). La diff\'erence avec la section 4 est que nous effectuons une chirurgie CR. Cela
ne pose toutefois pas de probl\`eme particulier car nous recollons une vari\'et\'e CR-triviale (voisinage de $\{0\}$ dans $\Bbb C^n$ 
produit $[0,\epsilon [$) le long d'un anneau produit.

\medskip

Cette op\'eration fabrique une famille de d\'eformations $\bar{\Cal W}$ de $\Bbb P^n$ param\'etr\'ee par l'intervalle. Comme $\Bbb P^n$ est rigide, cette famille est triviale au voisinage de la fibre centrale. En quittant la famille uniforme de
diviseurs et le voisinage de $0$, on en d\'eduit le th\'eor\`eme 2.
\medskip

Pour finir la preuve du th\'eor\`eme 3, on se sert de l'action. On poss\`ede en effet une uniformisation $\psi$ 
de $V_0$, compl\'ementaire d'un
voisinage de $0\in\bar {\Cal W}$, \`a valeurs 
dans $\Bbb C^n\setminus\{0\}\times [0,\infty[$. Sur le domaine d'uniformisation, on peut transporter
l'action de g\'en\'erateur $g$ en une action de g\'en\'erateur $h$; 
mais par Hartogs \`a param\`etres, cette action s'\'etend aux fibres compl\`etes de $\Bbb C^n\setminus\{0\}\times [0,\infty[$.
Remarquons en effet que le noyau de Cauchy que l'on peut utiliser ici pour faire l'extension est naturellement $C^\infty$ en le 
param\`etre transverse. 
\medskip
Ceci permet d'\'etendre $\psi$ \`a un domaine $V_1$ du type

\medskip
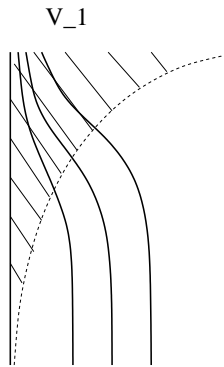
\hfil\scaledpicture 3.4in by 5.8in (domaineV1 scaled 300) \hfil
\medskip

en posant
$$
z\in V_1\longmapsto h^p \circ \psi\circ g^{-p}(z)\in\Bbb C^n\setminus \{0\}\times [0,\infty [
$$
o\`u $p$ est n'importe quel entier tel que $g^{-p}(z)$ appartienne \`a $V_0$.
$\square$

\head
{\bf 6. D\'eformations feuillet\'ees versus d\'eformations des vari\'et\'es de Hopf.}
\endhead

Le point de d\'epart de ce travail \'etait l'\'etude des d\'eformations des vari\'et\'es de Hopf. On sait \cite{Ha} que les vari\'et\'es 
de Hopf forment une classe stable par petites d\'eformations. La question centrale naturelle est

\proclaim{Question} Est-ce qu'une grande d\'eformation d'une vari\'et\'e de Hopf est une vari\'et\'e de Hopf?
\endproclaim

Nous param\'etrons ici nos grandes d\'eformations par l'intervalle $[0,1]$, mais la question est aussi signifiante en param\'etrant
par le disque complexe.

\medskip

Rappelons que la deuxi\`eme surface d'Hirzebruch s'obtient comme grande d\'e\-for\-mation de $\Bbb P^1\times\Bbb P^1$, bien que ce produit
soit rigide \cite{M-K}, et ce type de questions est en g\'en\'eral tr\`es difficile. Par ailleurs, on conna\^\i t des
structures complexes non-Hopf sur $\Bbb S^{2n-1}\times\Bbb S^1$.
 
\proclaim{D\'efinition}
Soit $X$ une vari\'et\'e complexe diff\'eomorphe \`a $\Bbb S^{2n-1}\times\Bbb S^1$. 
Nous dirons que $X$ est une vari\'et\'e de Brieskorn-Van de Ven de poids $a=(a_0,\hdots a_n)$ si son rev\^etement universel holomorphe
est l'hypersurface quasi-homog\`ene de $\Bbb C^{n+1}$ priv\'ee de $0$ d'\'equation
$$
W_a=\{z\in\Bbb C^{n+1}\setminus \{0\}\quad\vert\quad z_0^{a_0}+\hdots+z_n^{a_n}=0\}
$$
o\`u les $a_i$ sont des entiers naturels sup\'erieurs ou \'egaux \`a deux de PGCD valant un.
\endproclaim

On trouvera dans \cite{B-V} la construction de telles vari\'et\'es, en quotientant $W_a$ par une homoth\'etie pond\'er\'ee. 
Tous les choix de $a$ ne donnent pas forc\'ement une vari\'et\'e complexe diff\'eomorphe \`a $\Bbb S^{2n-1}\times\Bbb S^1$. On obtient
\'egalement des structures complexes sur des produits sph\`eres exotiques par cercle. C'\'etait d'ailleurs la motivation originale
de \cite{B-V}. En tout cas, \`a $n$ fix\'e, il existe une infinit\'e de $n$-uples $a$ donnant le produit standard.
\medskip
Remarquons que le rev\^etement universel d'une Hopf est 
$W_{(1,\hdots,1)}$, ce qui permet d'inclure le cas Hopf dans la construction de \cite{B-VdV} si l'on rel\^ache la condition que les $a_i$
doivent \^etre sup\'erieurs \`a deux. Une Brieskorn-Van de Ven n'est pas biholomorphe \`a une Hopf, car elles ont des rev\^etements 
universels holomorphiquement distincts. En fait,

\proclaim{Proposition}
Si le produit $(a_0-1)\cdots (a_n-1)$ est diff\'erent de $(b_0-1)\cdots (b_n-1)$, les vari\'et\'es abstraites $W_a$ et $W_b$ ne sont pas biholomorphes et m\^eme non biholomorphes dans aucun voisinage du bout $0$.
\endproclaim 

Nous remercions Dominique Cerveau qui nous a expliqu\'e l'argument suivant.

\demo{Preuve}
Montrons le r\'esultat par contrapos\'ee ; 
soit $f$ un biholomorphisme entre ces deux vari\'et\'es. On supposera que les $a_i$ sont tous de degr\'e
sup\'erieur ou \'egal \`a deux, tandis que les $b_i$ peuvent \^etre \'egaux \`a un pour inclure le cas Hopf. On supposera \'egalement
que $f$ est d\'efinie globalement, le cas local s'obtenant en rempla\c cant dans la suite de l'argument $\Bbb C^{n+1}$ par une boule 
de $\Bbb C^{n+1}$ centr\'ee en $0$ suffisamment petite.
\medskip
  
L'application $f$ s'\'etend en un biholomorphisme de $\bar W_a=W_a
\cup\{0\}$ dans $\bar W_b$ fixant $0$. Comme $\bar W_a$ est une hypersurface Stein de $\Bbb C^{n+1}$, la fl\`eche en cohomologie
$H^0(\Bbb C^{n+1},\Cal O)\to H^0(\Bar W_a,\Cal O)$ est surjective, autrement dit $f$ est la restriction d'une application $F$ de 
$\Bbb C^{n+1}$ dans lui-m\^eme ; de m\^eme $f^{-1}$ est la restriction d'une application $G$.
\medskip
Evidemment, $F$ et $G$ n'ont aucune raison d'\^etre des biholomorphismes de $\Bbb C^{n+1}$ et ils sont inverses l'un de l'autre uniquement
en restriction \`a $\bar W_a$ ou $\bar W_b$. N\'eanmoins, on peut \'ecrire {\it globalement}
$$
G\circ F\equiv Id+h
$$
o\`u $h$ est une \'equation de $\bar W_a$. Comme $h$ est de degr\'e au moins $2$, le th\'eor\`eme d'inversion locale implique que $F$ est
localement inversible. En particulier $F$ est un biholomorphisme d'un voisinage de $0$ dans $\Bbb C^{n+1}$ sur son image qui envoie 
la trace de $\bar W_a$ dans ce voisinage sur la trace de $\bar W_b$ dans l'image. Mais le th\'eor\`eme de fibration de Milnor \cite{Mi}
entra\^\i ne alors que les fibrations associ\'ees \`a $\bar W_a$ et $\bar W_b$ ont m\^eme nombre de Milnor, i.e. que 
le produit $(a_0-1)\cdots (a_n-1)$ est \'egal \`a $(b_0-1)\cdots (b_n-1)$.
$\square$
\enddemo

 Comme application des techniques pr\'ec\'edentes, on peut montrer :

\proclaim{Th\'eor\`eme 4}
Soit $\pi : \Cal X\to [0,1]$ une famille de d\'eformations dont les fibres $X_t$ sont des vari\'et\'es de Hopf pour tout $t$ diff\'erent de $1$. On suppose $n\geq 3$.
\medskip

Alors la fibre $X_1$ ne peut \^etre une vari\'et\'e de Brieskorn-Van de Ven. 
\medskip
De m\^eme, soit $\pi : \Cal X\to [0,1]$ une famille de d\'eformations dont les fibres $X_t$ sont des vari\'et\'es de Brieskorn de poids 
$a$ pour tout $t$ diff\'erent de $1$. Alors $X_1$ n'est ni une vari\'et\'e de Hopf, ni une Brieskorn-Van de Ven de nombre de Milnor 
distinct.
\endproclaim

A vrai dire, nous ne connaissons pas de vari\'et\'e complexe diff\'eomorphe \`a $\Bbb S^{2n-1}\times\Bbb S^1$ mis \`a part les vari\'et\'es de Hopf et les vari\'et\'es de Brieskorn-Van de Ven. Donc ce th\'eor\`eme laisse l'alternative suivante~: soit une grande
d\'eformation d'une Hopf est une Hopf, soit il existe une structure complexe ``exotique'' sur $\Bbb S^{2n-1}\times\Bbb S^1$ dans la
m\^eme classe de grande d\'eformation que les Hopf.

\demo{Preuve}
Pla\c cons-nous dans le cas Hopf, l'argument \'etant identique dans le deuxi\`eme cas.
Supposons le contraire. Comme le rev\^etement universel holomorphe $\tilde X_1$ de $X_1$ est une hypersurface quasi-homog\`ene $W_a$ de $\Bbb C^{n+1}$ priv\'ee de $0$, on peut conclure de la proposition 1 
que la d\'eformation $\tilde \Cal X$ de $\tilde X_1$ est localement pseudo-triviale. En effet, on a $\tilde X_1=\bar W_a\setminus\{0\}$ 
avec $0$ de codimension au moins trois dans $\bar W_a$ et de codimension homologique nulle (puisque $\bar W_a$ est
une intersection compl\`ete, cf \cite{A-G, Proposition 3}). On d\'eduit alors de \cite{Sc} l'isomorphisme entre $H^1(\tilde X_1,\Theta)$ et $H^1(\bar W_a, \Theta)$. Ce dernier groupe \'etant nul, puisque $\bar W_a$ est Stein, les hypoth\`eses
de la proposition 1 sont v\'erifi\'ees.
\medskip
Tout ceci entra\^{\i}ne par \cite{Ro} que, pour $t$ proche de $1$, les fibres $\tilde X_t$ peuvent \^etre modifi\'ees par
chirurgie au voisinage du bout $0$ en un espace de Stein $W_t$. Plus pr\'ecis\'ement, on construit $W_t$ en recollant \`a $\tilde X_t$ priv\'e d'un voisinage du bout $0$ un voisinage de $0$ dans $\bar W_a$. Notons en particulier
que tous les $W_t$ poss\`ede une unique singularit\'e et sont biholomorphes au voisinage de cette singularit\'e \`a un voisinage de $0$ dans $\bar W_a$.
\medskip

En raisonnant comme dans la preuve du th\'eor\`eme 1, comme $\tilde X_t$ supporte une contraction holomorphe, on a en fait une injection holomorphe de $\tilde X_t$ dans $W_t$ qui \'evite uniquement la singularit\'e. Ainsi
$\tilde X_t$ admet v\'eritablement une
compl\'etion Stein singuli\`ere par un point. Mais $\tilde X_t$ est isomorphe \`a $\Bbb C^n\setminus\{0\}$, donc admet l'\'evidente compl\'etion r\'eguli\`ere par un point $\Bbb C^n$ ; 
par unicit\'e d'une compl\'etion Stein, on aurait alors un biholomorphisme d'un voisinage de $0$ dans $\bar W_a$ \`a valeurs dans un 
voisinage de $0$ dans $\Bbb C^n$.
Contradiction avec la proposition pr\'ec\'edente.  
$\square$
\enddemo

Comparons maintenant avec le cas feuillet\'e. 
\medskip
\hfil\scaledpicture 9.2in by 3.4in (cylindre scaled 300) \hfil
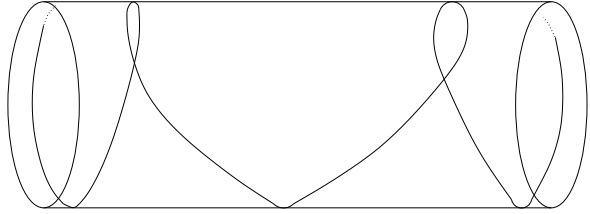
\medskip

\definition{D\'efinition}
Soit $M$ une vari\'et\'e lisse compacte de groupe fondamental $\Bbb Z$. Soient $X_0$ et $X_1$ deux vari\'et\'es compactes complexes diff\'eomorphes \`a $M$.
Une {\it d\'eformation feuillet\'ee} entre $X_0$ et $X_1$ est un feuilletage par vari\'et\'es complexes sur le cylindre $W=M\times [0,1]$ tel que :

\medskip
\noindent (i) Les deux composantes de bord sont les deux seules feuilles compactes ; l'une est biholomorphe \`a $X_0$ et l'autre \`a $X_1$.

\noindent (ii) Le rev\^etement universel feuillet\'e de $W$ est diff\'eomorphe \`a $\tilde M\times [0,1]$ feuillet\'e trivialement
par les niveaux de la projection sur $[0,1]$. De surcro\^\i t,  toutes les feuilles de l'int\'erieur de $W$ s'accumulent sur le bord (plus pr\'ecis\'ement un bout s'accumule sur une composante de bord, l'autre
bout sur l'autre). 

\noindent (iii) Les holonomies des composantes de bord sont plates, contractante d'un c\^ot\'e et dilatante de l'autre.
\enddefinition

Cela donne le dessin pr\'ec\'edent, en prenant pour $M$ un cercle.
On a 

\proclaim{Th\'eor\`eme 5}
Soit $X_1$ une vari\'et\'e compacte complexe obtenue par d\'eformation feuillet\'ee d'une vari\'et\'e de Hopf $X_0$. Si l'une des feuilles non compactes de la d\'e\-for\-ma\-ti\-on est isomorphe \`a $\Bbb C^n\setminus\{0\}$, alors 
\medskip
\noindent (i) $X_1$ est une vari\'et\'e de Hopf. 

\noindent (ii) Toutes les feuilles non compactes sont biholomorphes \`a $\Bbb C^n\setminus\{0\}$.

\noindent (iii) Le rev\^etement universel CR de la d\'eformation priv\'ee de $X_0$ et de $X_1$ est CR-isomorphe \`a $\Bbb C^n\setminus\{0\} \times ]0,1[$.



\endproclaim

\demo{Preuve}
Soit $(W, \Cal F)$ la d\'eformation feuillet\'ee joignant $X_0$ \`a $X_1$. D'apr\`es le point (ii) de la d\'efinition de d\'eformation feuillet\'ee, {\it toutes} les feuilles non compactes s'accumulent sur $X_0$ et sur $X_1$. C'est
donc en particulier le cas de la feuille biholomorphe \`a $\Bbb C^n\setminus\{0\}$. Cela suffit, d'apr\`es le th\'eor\`eme 1, pour montrer que $X_1$ est une vari\'et\'e de Hopf.
\medskip
Soit $\tilde W$ le rev\^etement universel de $W$. Notons qu'il est {\it diff\'eomorphe} \`a $\Bbb C^n\setminus\{0\} \times [0,1]$ et qu'il poss\`ede une projection $\pi$ sur $[0,1]$ dont les niveaux
d\'efinissent un feuilletage \`a feuilles complexes. D'apr\`es le corollaire 1, toutes les feuilles
non compactes sont de type $0$ (respectivement de type $\infty$) au voisinage de $X_1$. Et toujours d'apr\`es ce corollaire, toutes les feuilles non compactes sont de type $\infty$ (respectivement de type $0$) au voisinage de $X_0$.
On peut donc compactifier uniform\'ement $\tilde W$ priv\'e de $\tilde X_0$ et de $\tilde X_1$ en une famille de d\'eformations 
param\'etr\'ee par $]0,1[$ en utilisant deux fois le lemme de compactification uniforme. En utilisant le th\'eor\`eme de \cite{H-K}
comme dans la partie 5, on montre que toutes les feuilles sont biholomorphes \`a $\Bbb P^n$. 
Mais alors, comme {\it toutes} les fibres sont $\Bbb P^n$, la famille de d\'eformations construite est globalement triviale. La compactification effectu\'ee \'etant uniforme aux deux bouts, on en d\'eduit le point (iii) et donc 
automatiquement aussi le point (ii).
$\square$
\enddemo

Nous ne savons pas si la famille {\it \`a bord} de rev\^etements universels est triviale. En effet, chacune des deux compactifications
peut \^etre suppos\'ee uniforme d'un c\^ot\'e (i.e. en $0$ ou en $1$), mais pas de l'autre. On peut, en s'inspirant de ce que nous avons 
fait dans les paragraphes pr\'ec\'edents, compactifier $W$ en une famille de d\'eformations de vari\'et\'es compactes sur l'intervalle
ferm\'e $[0,1]$. Mais le ``d\'efaut'' d'uniformisation ne permet pas de conclure directement que cette famille est $\Bbb P^n\times [0,1]$.
Il faut penser que, partant de la famille triviale $\bar W=\Bbb P^n\times [0,1]$, si l'on enl\`eve une section $s:[0,1]\to W$ lisse
sur $]0,1[$ mais seulement continue en $0$ et en $1$, la famille $\bar W\setminus s$ n'est pas CR-triviale (cf \cite{M-V, \S 2}). 
 
\vfill
\eject

\enddocument
\bye